\documentclass{article}

\usepackage{graphicx}      
\usepackage{amsmath}
\usepackage{amssymb}
\usepackage{amsthm}
\usepackage{graphicx}
\usepackage{bm}
\usepackage{tikz}
\usepackage{xcolor}
\usepackage{authblk}
\usepackage{fullpage}

\newcommand{\grad}[0]{\textbf{grad}}
\renewcommand{\d}[0]{{\rm d}}

\newcommand{\Id}{\text{Id}}

\newtheorem{theorem}{Theorem}

\newtheorem{corollary}[theorem]{Corollary}

\newtheorem{remark}[theorem]{Remark}

\newenvironment{pf}[1][]{\begin{proof}#1\end{proof}}

\title{\bf On implicit and explicit representations for 1D distributed port-Hamiltonian systems\thanks{This work was partially supported  by the AID from the French
Ministry of the Armed Forces, and the IMPACTS project,
entitled IMplicit Port-HAmiltonian ConTrol Systems, funded by the French National Research Agency (project ANR-21-CE48-001),
{\tt https://impacts.ens2m.fr/}}} 
		
\author[1]{Antoine Bendimerad-Hohl\thanks{antoine.bendimerad-hohl@student.isae-supaero.fr}}
\author[1]{Denis Matignon\thanks{denis.matignon@isae.fr}}
\author[1]{Ghislain Haine\thanks{ghislain.haine@isae.fr}}
\author[2]{Laurent Lef{\`e}vre\thanks{laurent.lefevre@lcis.grenoble-inp.fr}}
		
\affil[1]{{\small Fédération ENAC ISAE-SUPAERO ONERA, Université de Toulouse, France}}
\affil[2]{{\small Univ. Grenoble Alpes, Grenoble INP, LCIS, 26000 Valence, France}}

\date{}

\begin{document}

\maketitle

\begin{abstract} \;
First, two examples of 1D distributed port-Hamiltonian systems with dissipation, given in explicit (descriptor) form, are considered: the Dzekster model for the seepage of underground water and a nanorod model with non-local viscous damping. Implicit representations in Stokes-Lagrange subspaces are formulated. These formulations lead to modified Hamiltonian functions with spatial differential operators. The associated power balance equations are derived, together with the new boundary ports. Second, the port-Hamiltonian formulations for the Timoshenko and the Euler-Bernoulli beams are recalled, the latter being a flow-constrained version of the former. Implicit representations of these models in Stokes-Lagrange subspaces and corresponding power balance equations are derived. Bijective transformations are proposed between the explicit and implicit representations. It is proven these transformations commute with the flow-constraint projection operator. 
\end{abstract}
		

\section{Introduction}
Port-Hamiltonian systems (pHs) have been developped for the modelling, (co-)simulation and control of
complex multiphysics systems (see \cite{van2014port} for an introductive textbook or \cite{duindam2009modeling} for applications of the approach to various domains, including mechanics, electromagnetism and irreversible thermodynamics). This framework has been extended to the case of distributed parameter systems (DPS) with boundary energy flows in the seminal paper \cite{van2002hamiltonian}. The semi-group approach, when applied to linear distributed parameter pHs, gives rise to particularly elegant and practical results in the 1D case (see \cite{jacob2012linear} for an introduction), it has been recently extended to $n$D systems in \cite{BHMCAM2023}. Since 2002, the literature on distributed pHs has grown considerably, with both theoretical and application papers (see \cite{rashad2020twenty} for a review).

Port-Hamiltonian system dynamics with algebraic constraints have been considered as well, leading to finite-dimensional port-Hamiltonian Differential-Algebraic Equations systems (pH-DAEs, see for instance \cite{beattie2018linear} for linear descriptor systems in matrix representation or \cite{mehrmann2023control} for a review on control applications). From a modeling or geometric perspective, these algebraic constraints arise either from interconnection of the subsystems composing the overall system, resulting in constraints between effort and flow variables in the underlying Dirac structure, or by implicit energy constitutive equations, resulting in constraints between effort and energy state variables assigned to live in some Lagrangian submanifold (see \cite{van2018generalized} or \cite{van2020dirac} for the nonlinear case). 

Recently, examples of distributed parameter models given in implicit form, either singular or not, have been considered  in the pHs settings (see for instance \cite{yaghi2022port} for an implicit formulation of the Allen-Cahn equation, \cite{jacob2022solvability} considering the Dzektser equation (seepage of underground water) or \cite{heidari2019port,heidari2022nonlocal} considering a nanorod with non-local visco-elastic constitutive equations), along with structure-preserving numerical methods, see \cite{bendimerad2023implicit}. \cite{maschke2023linear} extend boundary control pHs to a class of systems where the variational derivative of the Hamiltonian is replaced by a pair of reciprocal operators, generalizing -- in the infinite-dimensional settings -- the implicit definition of the energy by a so-called Stokes-Lagrange subspace associated to the reciprocal operators. This definition (of boundary pHs defined on a Stokes-Lagrange subspace) allows representing some of the previously cited distributed parameter models given in implicit form. Specifically, \cite{maschke2023linear} consider the example of an elastic rod with non-local elasticity relation.

In this paper, we propose first an application of this Stokes-Lagrange subspace approach to dissipative systems, in the case of local dissipation only (Dzekster equation, subsection \ref{subsec:_Dzektser}), and in the case with local and non-local dissipative ports (see the nanorod example in subsection \ref{subsec:_nanorod}). In both cases, classical explicit formulations and implicit representations using Stokes-Lagrange operators are proposed.  Note that in the considered implicit representations, spatial derivative operators are present inside the Hamiltonian as proposed for instance in \cite{scholberl2014auto}, where  infinite-dimensional pHs are defined on jet bundles. In \cite{preuster2024jet}, Boussinesq, elastic rod and Allen-Cahn equations are considered as examples of such systems. The authors propose a  lift in the jet space where the Hamiltonian density only depends on the extended state variable. Then, geometric formulations with Stokes-Dirac structures are applicable. 

In the second part of the paper, we derive similarly explicit and implicit representations for the Timoshenko and Euler-Bernoulli beams. The Euler-Bernoulli model may be seen as flow-constrained version of the explicit Timoshenko beam and is therefore an example of constrained pHs with an implicit definition of the energy via a Stokes-Lagrange subspace. Bijective operators which transform explicit formulations to implicit ones are proposed for the Timoshenko beam model and for the Euler-Bernouilli reduced model. These operators are similar to transformations between DAEs and geometric representations analyzed in \cite{mehrmann2023differential}, in the finite-dimensional LTI case. It is shown that the reduction (from Timoshenko to Euler-Bernoulli models) and the bijective explicit-implicit transformations commute in this particular example. This allows, for instance, to consider the implicit formulation for mathematical analysis, while the constrained explicit formulation (DAEs) is used for numerical computations.

\section{Preliminary results}

\subsection{Variational derivatives}

\begin{theorem}

Let   $\mathcal{K}:D(\mathcal{K}) \subseteq L^2(\Omega,\mathbb{R}^n) \rightarrow L^2(\Omega,\mathbb{R}^m) $ be a closed and densely-defined (differential) linear operator, and let
\begin{equation*}
   \forall \alpha \in D(\mathcal{K}), \quad \mathcal{H}(\alpha) := \frac{1}{2}\int_\Omega |\mathcal{K}(\alpha)|^2 \d x,
\end{equation*}
be an energy functional.

Assume that the following abstract Green's identity holds for all $\alpha \in D(\mathcal{K}), \, \beta \in D(\mathcal{K}^\dag)$
\begin{equation}
\label{eq:Green-identity}
\int_\Omega \mathcal{K} (\alpha) \cdot \beta 
=
\int_\Omega \alpha \cdot \mathcal{K}^\dag (\beta)
+ \left\langle \gamma(\alpha), \mathcal{C}(\beta) \right\rangle_{\mathcal{U}, \mathcal{Y}},
\end{equation}
where $\mathcal{K}^\dag:D(\mathcal{K}^\dag) \subseteq L^2(\Omega,\mathbb{R}^m) \rightarrow L^2(\Omega,\mathbb{R}^n)$ is another closed and densely-defined linear operator, called the \emph{formal adjoint} of $\mathcal{K}$, $\gamma \in \mathcal{L}(D(\mathcal{K}), \mathcal{U})$ is a (boundary) control operator on the Hilbert space~$\mathcal{U}$, and $\mathcal{C} \in \mathcal{L}(D(\mathcal{K}^\dag), \mathcal{Y})$ its colocated observation operator, with $\mathcal{Y} := \mathcal{U}'$.

Then the \emph{weak} variational derivative of $\mathcal{H}$ with respect to $\alpha$ exists and is given as
\begin{equation*}
    \delta^w_{\alpha} \mathcal{H} (\alpha) = \mathcal{K}^\dag(\mathcal{K} (\alpha)).
\end{equation*}
\end{theorem}
\begin{pf}
    Let $\alpha \in D(\mathcal{K}), \varepsilon \in \mathbb{R},$ then, for all $\beta \in C^\infty_c(\Omega, \mathbb{R}^m) =: D(\Omega)$, one has
    \begin{equation*}
        \begin{aligned}
            \mathcal{H}(\alpha \! + \! \varepsilon \beta ) =& \frac{1}{2} \int_\Omega \mathcal{K}(\alpha + \varepsilon \beta ) \cdot \mathcal{K}(\alpha + \varepsilon \beta) \d x \\
            =& \frac{1}{2} \int_\Omega \big( | \mathcal{K}(\alpha) |^2 + 2\varepsilon \, \mathcal{K}(\alpha) \cdot  \mathcal{K}(\beta) \\ 
            & \quad \quad + \varepsilon^2 \, | \mathcal{K}(\beta) |^2 \big) \d x \\
            =& \mathcal{H}(\alpha)  + \varepsilon \int_\Omega \alpha \cdot \mathcal{K}^\dag(\mathcal{K}(\beta)) \d x  + o(\varepsilon).
        \end{aligned}
    \end{equation*}
Defining the weak variational derivative in the sense of distribution (thanks to~\eqref{eq:Green-identity})
$$
\left\langle \delta^w_{\alpha} \mathcal{H}(\alpha), \beta \right\rangle_{D(\Omega)', D(\Omega)} := \int_\Omega \alpha \cdot \mathcal{K}^\dag(\mathcal{K}(\beta)) \d x,
$$
gives the result.
\end{pf}

\begin{corollary} \label{corollary-variational-derivative}
Given the Hamiltonians $\mathcal{H}_1$ and $\mathcal{H}_2$, defined on $H^1(\Omega, \mathbb{R}), H^2(\Omega, \mathbb{R})$, respectively, by
$$
    \mathcal{H}_1(\alpha) = \frac{1}{2} \int_\Omega (\partial_x \alpha)^2 \d x, \, 
    \mathcal{H}_2(\alpha) = \frac{1}{2} \int_\Omega (\partial_{x^2}^2 \alpha)^2 \d x,
$$
their weak variational derivatives are given by
$$ 
\begin{aligned}
    \delta^w_{\alpha} H_1(\alpha) = - \partial_{x^2}^2 \alpha, \quad
    \delta^w_{\alpha} H_2(\alpha) = \partial_{x^4}^4 \alpha, 
\end{aligned}
$$
respectively.
\end{corollary}

\begin{pf}
Take $\mathcal{K} = \partial_x$, $\mathcal{K}^\dag = - \partial_x$ in the first case, and $\mathcal{K} = \mathcal{K}^\dag = \partial^2_{x^2}$ in the second case.
\end{pf}

\subsection{Operator transposition}

In the following, the superscript notations on a variable $x^E,\,x^I$ refers to the explicit or implicit representation of a system, respectively.

Let us consider an operator $\mathcal{K}$ as before, a set of distributed second order tensors $\bm{\kappa} \in L^\infty(\Omega,M_n(\mathbb{R})), \, \bm{\eta} \in L^\infty(\Omega, M_m(\mathbb{R}))$, with $\forall x \in \Omega, \, \bm{\kappa}(x) = \bm{\kappa}^\top(x)  \geq \kappa > 0, \bm{\eta}(x) = \bm{\eta}^\top(x) \geq \eta > 0$ almost everywhere, and a pHs defined as
\begin{equation*}
        \begin{bmatrix}
            \partial_t \alpha^E_1 \\
            \partial_t \alpha^E_2
        \end{bmatrix} = \begin{bmatrix}
            0 & \mathcal{K}^\dag \\ - \mathcal{K} & 0
        \end{bmatrix} \begin{bmatrix}
            e_1 \\ e_2
        \end{bmatrix}, \quad
    \begin{cases}
        e^E_1 = \delta_{\alpha_1} \mathcal{H}^E, \\
        e^E_2 = \delta_{\alpha_2} \mathcal{H}^E,
    \end{cases}
\end{equation*}
with the corresponding Hamiltonian defined as
$$ 
\mathcal{H}^E = \frac{1}{2} \int_\Omega \alpha^E_1 \cdot \bm{\kappa}\, \alpha^E_1 +\alpha^E_2 \cdot \bm{\eta}\, \alpha^E_2 \d x. 
$$

Let us now consider a second port-Hamiltonian system defined as
\begin{equation*}
        \begin{bmatrix}
            \partial_t \alpha^I_1 \\
            \partial_t \alpha^I_2
        \end{bmatrix} = \begin{bmatrix}
            0 & \Id \\ - \Id & 0
        \end{bmatrix} \begin{bmatrix}
            e_1^I \\ e_2^I
        \end{bmatrix}, \quad
    \begin{cases}
        e^I_1 = \delta_{\alpha^I_1} \mathcal{H}^I, \\
        e^I_2 = \delta_{\alpha^I_2}^w \mathcal{H}^I,
    \end{cases}
\end{equation*}
with the corresponding Hamiltonian defined as
$$
\mathcal{H}^I = \frac{1}{2} \int_\Omega \alpha^I_1 \cdot \bm{\kappa}\, \alpha^I_1 + \mathcal{K}(\alpha^I_2) \cdot \bm{\eta}\, \mathcal{K}(\alpha^I_2) \d x.
$$
For the sake of readability, weak variational derivative will be denoted with the same symbol as usual variational derivative $\delta$ in the sequel of the paper.

Finally, let us define an operator $\mathcal{G}:  L^2(\Omega,\mathbb{R}^n) \times D(\mathcal{K}) \rightarrow L^2(\Omega,\mathbb{R}^n)\times L^2(\Omega,\mathbb{R}^m) $ and its $\dag$-companion
$$
\mathcal{G} := \begin{bmatrix}
    \Id & 0 \\ 0 & \mathcal{K}
\end{bmatrix},
\quad
\text{and}
\quad
\mathcal{G}^\dag := \begin{bmatrix}
    \Id & 0 \\ 0 & \mathcal{K}^\dag
\end{bmatrix},
$$
and the structure matrices of operators
$$ \mathcal{J}^E = \begin{bmatrix}
    0 & \mathcal{K}^\dag \\ - \mathcal{K} & 0
\end{bmatrix}, 
\quad
\mathcal{J}^I = \begin{bmatrix}
    0 & \Id \\ - \Id & 0
\end{bmatrix}. $$

We then get the following theorem:
\begin{theorem} The previously defined operator $\mathcal{G}$ allows passing from one representation to the other, the transformation being given as
    \begin{equation*}
        \begin{cases}
        \mathcal{G} \alpha^I =  \alpha^E, \\
        e^I = \mathcal{G}^\dag e^E, \\
        \mathcal{J}^E = \mathcal{G}\mathcal{J}^I\mathcal{G}^\dag.
        \end{cases}
    \end{equation*}
    The two systems are said to be \emph{equivalent}.
\end{theorem}
\begin{pf}
    A direct computation gives the results.
\end{pf}

\begin{remark}
    Note that given such a transformation $\mathcal{G}$, each state $\alpha^I \in \ker(\mathcal{G})$ does not participate in the Hamiltonian, in particular they are removed by the transformation. This explains why when writing the wave equation in pH formulation, the deformation $\bm{\varepsilon}$ is used instead of the displacement $w$:  the Hamiltonian $\mathcal{H} = \frac{1}{2}\int_\Omega p^2 + |\grad(w)|^2 \d x$ only depends on the momentum $p$ and deformation $\grad(w) = \bm{\varepsilon}$, hence the transformation removing $\grad$ from the Hamiltonian yields $(p,\bm{\varepsilon})$ as the state. And because the kernel of $\grad$ is the set of constant functions, rigid body motion is lost during the transformation.
\end{remark}

\section{Two examples with damping}
In \S~\ref{subsec:_Dzektser}, the seepage model is considered, and in \S~\ref{subsec:_nanorod}, the example of the nanorod is studied.
Throughout this paper, the 1D~domain $\Omega$ is defined as a bounded interval $\Omega = [a,b]$.

\subsection{Dzektser}
\label{subsec:_Dzektser}
Following \cite{Dzektser72}, let us consider the following seepage model of underground water in 1D
$$
(\Id - \varepsilon^2 \partial_{x^2}^2)\,\partial_t h = \alpha\,\partial_{x^2}^2 h - \beta\,\partial_{x^4}^4 h,
$$
with $\alpha >0, \beta>0$. The system admits a port-Hamiltonian representation given as
\begin{equation}\label{eqn:dzektser-explicit}
    \begin{cases}
        \begin{bmatrix}
            (\Id - \varepsilon^2 \partial_{x^2}^2)  & 0 & 0\\
            0 & \Id & 0\\
            0 & 0 & \Id
        \end{bmatrix}\!\! \begin{bmatrix}
            h \\
            F_\nabla \\
            F_\Delta
        \end{bmatrix} \!\! =\!\! \begin{bmatrix}
            0 & \partial_x & - \partial_{x^2}^2 \\
            \partial_x & 0 & 0 \\
            \partial_{x^2}^2 & 0 & 0
        \end{bmatrix}\!\!\!\begin{bmatrix}
            h \\
            E_\nabla \\
            E_\Delta
        \end{bmatrix}\!\!,\\
        \text{with the resistive relations} \\
        E_\nabla = \alpha F_\nabla, \\
        E_\Delta = \beta F_\Delta.
    \end{cases}
\end{equation}
Following \cite{maschke2023linear} one can then define the Lagrange subspace operators
\begin{equation*}
    \mathcal{S} = (\Id - \varepsilon^2 \partial_{x^2}^2), \quad \mathcal{P} = \Id,
\end{equation*}
the Dirac structure operator
$$ 
\mathcal{J} = \begin{bmatrix}
            0 & \partial_x & - \partial_{x^2}^2 \\
            \partial_x & 0 & 0 \\
            \partial_{x^2}^2 & 0 & 0
        \end{bmatrix}, 
$$
and the resistive structure operator: $\mathcal{R} = \begin{bmatrix}
\alpha & 0 \\ 0 & \beta    
\end{bmatrix}$.
Finally, following \cite[Section 5.]{maschke2023linear}, the Hamiltonian is given as 
$$ \begin{aligned}
\mathcal{H} &= \frac{1}{2}\int_\Omega h \mathcal{S}^\dag \mathcal{P} h \d x + \frac{1}{2}[\varepsilon^2 h \partial_x h]_a^b, \\
            &= \frac{1}{2} \int_\Omega h^2 + \varepsilon^2\,(\partial_x h)^2 \d x. 
\end{aligned}
$$
\begin{theorem} (Dzektser Power balance) The power balance related to~\eqref{eqn:dzektser-explicit} reads:
\begin{equation*}\begin{aligned}
    \frac{\d}{\d t} \mathcal{H} = & \,\, [\alpha h  \partial_x h - \beta  h \,  \partial_{x^3}^3 h + \beta \partial_x h \, \partial_{x^2}^2 h ]_a^b \\
                                   - &\int_\Omega \!\!\! \left (a (\partial_x h)^2 + \beta(\partial_{x^2}^2h)^2 \right ) \d x + [\varepsilon^2 \partial_x h \,  \partial_t h]_a^b.
\end{aligned}
\end{equation*}    
\end{theorem}

Making use of  $\rm{tr}$ as  the Dirichlet trace operator, one can then define
 $(f_\partial,e_\partial)$ the \emph{\textbf{power} boundary port} related to the Stokes-Dirac structure as
\begin{equation*}
    \begin{aligned}
        f_\partial    =& \text{tr}(\begin{bmatrix}\,
            h, & &\, h, \,& &\partial_x h \,
        \end{bmatrix}),  \\
        e_\partial =& \text{tr}(\begin{bmatrix}\,
            \alpha \partial_x h,  && -\beta \partial_{x^3}^3 h, && \beta \partial_{x^2}^2 h \,
        \end{bmatrix}),
    \end{aligned}
\end{equation*}
 and $(\chi_\partial,\varepsilon_\partial)$ the  \emph{\textbf{energy} boundary port} related to the Stokes-Lagrange subspace as
\begin{equation*}
    \begin{aligned}
        \chi_\partial =\text{tr}(\begin{matrix}
            h
        \end{matrix}), \qquad
        \varepsilon_\partial =\text{tr}(\begin{matrix}
            \varepsilon^2 \partial_x h
        \end{matrix}).
    \end{aligned}
\end{equation*}
The power balance then reads
\begin{eqnarray*}
    \frac{\d}{\d t} \mathcal{H} &=& [f_\partial \cdot e_\partial + \frac{\d }{\d t}(\chi_\partial) \varepsilon_\partial]_a^b - \int_\Omega
        (\alpha F_\nabla^2 + \beta F_\Delta^2)\, \d x\,,\\
        &\leq&  [f_\partial \cdot e_\partial + \frac{\d }{\d t}(\chi_\partial) \varepsilon_\partial]_a^b \,.
\end{eqnarray*}
\begin{remark}
    The previous inequality is the extension to the case of lossy systems  of the equality established in \cite{maschke2023linear} for the case of lossless systems.
\end{remark}
\begin{remark}
The energy boundary port $(\chi_\partial, \varepsilon_\partial)$ vanishes when $\varepsilon \rightarrow 0$, i.e., when the nonlocal term is removed, leading to a classical dissipative pHs.
\end{remark}
\subsection{Nanorod}
\label{subsec:_nanorod}
Let us first begin by  writing  both versions of the Nanorod example and then compare them.	
	\subsubsection{Explicit representation}
Following \cite{heidari2019port}, the Hamiltonian of the system reads
\begin{equation*}
 \begin{aligned}
		\mathcal{H} 
		:= \frac{1}{2} \int_\Omega &
		a^2 w^2 
		+ \rho A \left( \partial_t w \right)^2 
		+ \mu \rho A \left( \partial^2_{tx}w \right)^2 
		\\& \!\!\!+ \left( E A + \mu a^2 \right) \left(\partial_x w \right)^2 \d x,
  \end{aligned}
	\end{equation*}
and the state variable is given as
$$
z := 
 [
	w, \,
	\rho A \partial_t w, \,
	\mu \rho A \partial^2_{tx} w, \,
	\partial_x w, \,
	N
     ]^\top,
$$ 
with $\omega$ is the displacement, $\rho A \partial_t w$ the momentum density,  $\mu \rho A \partial^2_{tx}$ the flow variable of the non locality, $\partial_x w$ the strain and $N$ the stress resultant. Let us now define ${\cal E} := \text{Diag}(\Id, \Id, \Id, \Id, 0)$ and
$$
	{\cal Q} := 
	\begin{bmatrix}
		a^2 & 0 & 0 & 0 & 0 \\
		0 & \frac{1}{\rho A} & 0 & 0 & 0 \\
		0 & 0 & \frac{1}{\mu \rho A} & 0 & 0 \\
		0 & 0 & 0 & EA + \mu a^2 & 0 \\
		0 & 0 & 0 & 0 & \Id 
	\end{bmatrix},
$$
which allows us to rewrite the Hamiltonian $\mathcal{H}$ as $\mathcal{H} = \frac{1}{2} \int_\Omega z^\top {\cal E}^\top {\cal Q} z$, with the algebraic property 
$$ 
{\cal E}^\top {\cal Q}  = {\cal Q}^\top {\cal E}.
$$
Defining furthermore
$$\begin{aligned}
	{\cal J} &:=
	\begin{bmatrix}
		0 & \Id & 0 & 0 & 0 \\
		-\Id & 0 & 0 & 0 & \partial_x \\
		0 & 0 & 0 & -\Id & \Id \\
		0 & 0 & \Id & 0 & 0 \\
		0 & \partial_x & -\Id & 0 & 0 
	\end{bmatrix},\\ 
	{\cal R} &:= \text{Diag}(
	    0,  b^2, \tau_d EA + \mu b^2, 0, 0),
 \end{aligned}
$$
where $b$ is the damping coefficient of the viscoelastic layer and $\tau_d$ the viscous damping of the nanorod.

The dynamics of the system is given by
\begin{equation}\label{eq:dynamic}
	{\cal E} \partial_t z = ({\cal J}-{\cal R})\,e, \qquad e = {\cal Q} z.
\end{equation}

\subsubsection{Implicit representation}

Let us now write  this system as an implicit port-Hamiltonian system, i.e., as a state-space representation where differential operators are present inside the Hamiltonian
\begin{equation}
	\begin{bmatrix}
		\partial_t w \\
		\partial_t \varepsilon \\
		\partial_t p \\
		f_d \\
		f_\sigma
	\end{bmatrix} = \begin{bmatrix}
	0 & 0 & \Id & 0 & 0\\
	0 & 0 & \partial_x & 0 & 0\\
	-\Id & \partial_x  & 0 & -\Id & \partial_x \\
	0 & 0 & \Id & 0 & 0 \\
	0 & 0 & \partial_x & 0 & 0
\end{bmatrix} \begin{bmatrix}
F \\  \sigma \\v \\ e_d \\ e_\sigma
\end{bmatrix},
\end{equation}
with $p$ the momentum, $w$ the displacement, $\varepsilon$ the strain, $f_d$ and $e_d$ the local and nonlocal dissipative ports, $v$ the velocity, $F$ the force of the media applied to  the nanorod, $\sigma$ the nonlocal stress, and $e_d,e_\sigma$ the local and nonlocal efforts linked to the dissipative port. The constitutive relations are
\begin{equation*}
	\begin{aligned}
		F = a w, &\quad  
		(\Id- \mu\,\partial_{x^2}^2) \sigma = E \varepsilon, \\
		v = \frac{p}{\rho A}, \quad
		&e_d = b f_d, \quad
		(\Id - \mu\,\partial_{x^2}^2) e_\sigma = \tau_d f_\sigma.
	\end{aligned}
\end{equation*}

Notice that $(\Id-\mu\,\partial_{x^2}^2)$ is now inside the constitutive relations, and that two constitutive relations have become non-local.

Let us denote by $z = (w,\varepsilon,p)$ the state and $e_S = (F,\sigma,v)$ the effort variable corresponding to the storage port.

 We can now define the $\mathcal{S}$ and $\mathcal{P}$ matrices as
\begin{equation}
		\mathcal{P}  := \begin{bmatrix}
		\Id & 0 & 0 \\
		0 & (\Id- \mu \,\partial_{x^2}^2) &0   \\
		0 & 0 & \Id \\
	\end{bmatrix}, \quad \mathcal{S}  := \begin{bmatrix}
	a & 0 & 0  \\
	0  & E & 0  \\
	0 & 0 & \frac{1}{\rho} \\
\end{bmatrix}.
\end{equation}
And we get the following relation between state and effort variables
$$ 
\mathcal{S}^\dag z = \mathcal{P}^\dag e_S. 
$$

In particular, we have that $\mathcal{S}^\dag \mathcal{P} = \mathcal{P}^\dag \mathcal{S}$, which shows that the corresponding Hamiltonian is non-negative. 
Moreover, we can define the resistive matrices
$$
\mathcal{R}_L = \begin{bmatrix}
	\Id & 0 \\
	0 & (\Id - \mu \,\partial_{x^2}^2)
\end{bmatrix},  \qquad \mathcal{R}_R = \begin{bmatrix}
b & 0 \\0  & \tau_d
\end{bmatrix},
$$
which yields the following \textbf{implicit} resistive structure
$$ 
\mathcal{R}_L^\dag \mathcal{R}_R  \geq 0, \mbox{ and } \mathcal{R}_L e_r = \mathcal{R}_R f_r,
$$
with $f_r = (f_d,f_\sigma)^\top$ and $e_r = (e_d, e_\sigma)^\top$.

Let us finally define the structure matrix
$$ 
\mathcal{J} = \begin{bmatrix}
	0 & 0 & \Id & 0 & 0\\
	0 & 0 & \partial_x & 0 & 0\\
	-\Id & \partial_x  & 0 & -\Id & \partial_x \\
	0 & 0 & \Id & 0 & 0 \\
	0 & 0 & \partial_x & 0 & 0
\end{bmatrix},
$$
then, the system becomes
$$ 
	\begin{bmatrix}
	\partial_t z \\
	f_r
\end{bmatrix} = \mathcal{J} \begin{bmatrix}
e_S \\ e_r
\end{bmatrix}, \quad
\begin{cases}
   \mathcal{R}_L e_r = \mathcal{R}_R f_r , \\
   \mathcal{S}^\dag z = \mathcal{P}^\dag e_S .
\end{cases} 
$$

Let us now write it using the \emph{image representation} of the Lagrange subspace
$$ 
	\begin{bmatrix}
		\mathcal{P} \partial_t \xi \\
		f_r
	\end{bmatrix} = \mathcal{J} \begin{bmatrix}
		e_S \\ e_r
	\end{bmatrix}, \quad
\begin{cases}
	\mathcal{R}_L e_r = \mathcal{R}_R f_r , \\
	e_S= 	\mathcal{S}  \xi .
\end{cases}
$$
with $\mathcal{P}\xi = z$ being the latent space variable. 
Following \cite[Section 5.]{maschke2023linear}, this representation allows us to define the corresponding Hamiltonian
$$
\begin{aligned}
\mathcal{H}^I := & \frac{1}{2} \int_\Omega \xi^\top \mathcal{P}^\dag \mathcal{S} \xi \d x + \frac{1}{2} [\mu E \xi_2  \partial_x \xi_2]_a^b, \\
                = & \frac{1}{2} \int_\Omega a{\xi_1}^2 + E {\xi_2}^2 + \mu E (\partial_{x} \xi_2)^2 + \frac{1}{\rho} {\xi_3}^2\,\d x.
\end{aligned}
$$
 
\begin{theorem}(Nanorod) The power balance reads
    \begin{eqnarray}
        \frac{\d}{\d t} \mathcal{H^I} &=&\,\, [\,\mu E \partial_t \xi_2 \partial_x \xi_2 + \sigma v + e_\sigma v - \mu e_\sigma \partial_x e_\sigma\,]_a^b \nonumber \\
        && - \int_\Omega b f_d^2 + \frac{1}{\tau_d} (e_\sigma^2 + \mu (\partial_x e_\sigma)^2) \d x.
    \end{eqnarray}
\end{theorem}
One can then identify the \emph{\textbf{power} boundary ports} $(f_\partial,e_\partial)$ and \emph{\textbf{energy} boundary ports} $(\chi_\partial,\varepsilon_\partial)$ as
\begin{equation*}
    \begin{aligned}
        f_\partial = \text{tr}(\begin{bmatrix}
            v, & v, & - \mu \frac{1}{\tau_d} \partial_x e_\sigma  
        \end{bmatrix}), \; & e_\partial = \text{tr}(\begin{bmatrix}
            \sigma, & e_\sigma, & e_\sigma
        \end{bmatrix}), \\
        \chi_\partial = \text{tr} (\begin{matrix}
            \xi_2
        \end{matrix}), \; & \varepsilon_\partial = \text{tr}(\begin{matrix}
            \mu E \partial_x \xi_2
        \end{matrix}). 
    \end{aligned}
\end{equation*}
The power balance then reads
\begin{eqnarray*}
     \frac{\d}{\d t} \mathcal{H^I} &= & \, \,   [\, f_\partial \cdot e_\partial + \frac{\d}{\d t}(\chi_\partial) \varepsilon_\partial\,]_a^b \\ 
    &&- \int_\Omega b f_d^2 + \frac{1}{\tau_d} (e_\sigma^2 + \mu (\partial_x e_\sigma)^2) \d x. \\
    &\leq &  \, \, [\, f_\partial \cdot e_\partial + \frac{\d}{\d t}(\chi_\partial) \varepsilon_\partial\,]_a^b\,.
\end{eqnarray*}

\section{Equivalent representations for classical beam models}

In \S~\ref{ss-ExpRep}, explicit representations both for the Timoshenko and the Euler-Bernoulli beam are recalled, then 
in \S~\ref{ss-ImpRep}, implicit representations are derived, thus recovering the examples treated in the jet bundle formalism first presented in \cite{scholberl2014auto} for Timoshenko, and in \cite[Example 5.1]{SchoberlSchlacherMathMod2015} for Euler-Bernoulli.
Finally in \S~\ref{ss-EquivRep}, the equivalence between these representations is proved, and the passage to the limit from Timoshenko to Euler-Bernoulli well understood in a common setting.

\subsection{Explicit representations}
\label{ss-ExpRep}
 \subsubsection{Timoshenko}

	The Timoshenko beam equation reads \cite{ducceschi2019conservative}	
$$
 \begin{cases}
		\rho A \partial_t^2 w = T_0 \partial_{x^2}^2 w + A \kappa G \partial_x (\partial_x w - \phi), \\
		\rho I \partial_t^2 \phi = EI \partial_{x^2}^2 \phi + A \kappa G (\partial_xw - \phi),
	\end{cases}
 $$
	with $\rho$ the mass density, $A$ the cross section area, $T_0$ the tension, $\kappa$ the sheer coefficient, $G$ the shear modulus, $I$ moment of inertia,  $E$ the young's modulus,   $w$ the transverse displacement and  $\phi$ the shear angle. Let us firstly write it as an explicit port-Hamiltonian system
    \begin{equation} \label{eqn:dynamic-timoshenko-explicit} \begin{bmatrix}
		\partial_t \varepsilon_w  \\
		\partial_t p_w \\
		\partial_t \varepsilon_\phi \\
		\partial_t p_\phi \\
		\partial_t \varepsilon_{w,\phi}
	\end{bmatrix} = \underbrace{\begin{bmatrix}
	0 & \partial_x & 0 &0& 0 \\
	\partial_x & 0 & 0 & 0  &\partial_x \\
	0 & 0 & 0  & \partial_x & 0\\
	0 & 0 & \partial_x & 0  & \Id \\
	0 & \partial_x & 0 & -\Id & 0 
\end{bmatrix}}_{=:\mathcal{J}^E} \begin{bmatrix}
	\sigma_w \\
	v \\
	\sigma_\phi \\
	\omega \\
	N
\end{bmatrix},
\end{equation}
with $p_w$ the linear momentum, $p_\phi$ the angular momentum, $\varepsilon_w = \partial_x w$ the deformation, $\varepsilon_\phi = \partial_x \phi$ the spatial derivative of the shear angle, and $\varepsilon_{w,\phi} = \partial_x w - \phi$ the difference between the deformation and the shear angle.

With the following constitutive relations
$$ 
\begin{aligned}
    v = \frac{p_w}{\rho A}, \quad   \sigma_w = T_0 \varepsilon_w, \quad \omega = \frac{p_\phi}{\rho I} \\
    \sigma_\phi = E I \varepsilon_\phi, \quad N = A \kappa G \varepsilon_{w,\phi},
\end{aligned}
$$
 and Hamiltonian
 \begin{equation}\label{eqn:ham-timoshenko-explicit}
    \begin{aligned}
     \mathcal{H}^E := \frac{1}{2} \int_\Omega& \Big( \frac{p_w}{\rho A}^2 + T_0 \varepsilon_w^2 \\& + \frac{p_\phi}{\rho I}^2 + E I \varepsilon_\phi^2 + A \kappa G \varepsilon_{w,\phi}^2 \Big) \d x.
    \end{aligned}
 \end{equation}

 Moreover, one can define the matrices 
 $$
 \mathcal{S}^E = \Id, \quad \mathcal{P}^E = \text{Diag}\begin{pmatrix}
     T_0, & \frac{1}{\rho A}, & EI, & \frac{1}{\rho I}, & A \kappa G
 \end{pmatrix},
 $$
with $\mathcal{P}^{E\top} \mathcal{S}^E = \mathcal{S}^{E\top} \mathcal{P}^E > 0.$ Finally one gets the following power balance:

 \begin{theorem}(Timoshenko power balance - explicit case) The power balance of \eqref{eqn:dynamic-timoshenko-explicit} reads
 \begin{equation}
     \frac{\d }{\d t} \mathcal{H}^E = \left [v \sigma_w + vN + \omega \sigma_\phi \right ]_a^b.
 \end{equation}
 One can identify the \emph{\textbf{power} boundary port} variables $f_\partial = \text{tr}(\begin{bmatrix}
     v & v & \omega
 \end{bmatrix}^\top), e_\partial = \text{tr}(\begin{bmatrix}
     \sigma_w & N &\sigma_\phi
 \end{bmatrix}^\top)$, tr being the trace operator; yielding $\frac{\d}{\d t}\mathcal{H}^E = [f_\partial \cdot e_\partial]_a^b $.
 \end{theorem}

\subsubsection{Euler--Bernoulli}

In order to reduce the Timoshenko model to the Euler--Bernoulli model, one simply needs to transform the ODE into a DAE by setting $\partial_t p_\phi$ and $\partial_t \varepsilon_{w,\phi}$ to zero in~\eqref{eqn:dynamic-timoshenko-explicit}
\begin{equation} \label{eqn:dynamic-eulerbernoulli-explicit-DAE}
\begin{bmatrix}
		\partial_t \varepsilon_w  \\
		\partial_t p_w \\
		\partial_t \varepsilon_\phi \\
		\bm{0} \\
		\bm{0}
	\end{bmatrix} = \mathcal{J}^E \begin{bmatrix}
	\sigma_w \\
	v \\
	\sigma_\phi \\
	\omega \\
	N
\end{bmatrix}. 
\end{equation}

 One can then define the constrained Hamiltonian
\begin{equation}\label{eqn:ham-eulerbernoulli-explicit}\begin{aligned}
    \mathcal{H}^E_c = \frac{1}{2} \int_\Omega & \frac{p_w}{\rho A}^2 + T_0 \varepsilon_w^2  + E I \varepsilon_\phi^2 \d x, 
    \end{aligned}
\end{equation}
and get the following power balance
\begin{theorem}(Euler-Bernoulli power balance - explicit case) The power balance of \eqref{eqn:dynamic-eulerbernoulli-explicit-DAE} reads
\begin{equation}
    \frac{\d }{\d t} \mathcal{H}^E_c = \left[ v \sigma_w - v \partial_x \sigma_\phi + \sigma_\phi \partial_x v \right]_a^b.
\end{equation}
\end{theorem}
One can identify the \emph{\textbf{power} boundary port} variables $f^c_\partial = \text{tr}(\begin{bmatrix}
     v & v &  \partial_x v
 \end{bmatrix}^\top\!), e^c_\partial \!= \!\text{tr}(\begin{bmatrix}
     \sigma_w & - \partial_x \sigma_\phi &\sigma_\phi
 \end{bmatrix}^\top \!)$.   This yields $\frac{\d}{\d t}\mathcal{H}^E = [f_\partial^c \cdot e_\partial^c]_a^b $.
 
\begin{remark}
One can solve the constraints analytically:
    $ \partial_x \sigma_\phi  = - N, \quad \partial_x v = \omega. $
Which yields the following reduced system
\begin{equation} \label{eqn:dynamic-eulerbernoulli-explicit}
    \begin{bmatrix}
         \partial_t \varepsilon_w \\ \partial_t \varepsilon_\phi \\ \partial_t p_w
    \end{bmatrix} = \underbrace{\begin{bmatrix}
        0 & 0 & \partial_x \\
        0  &0  & \partial_{x^2}^2 \\
         \partial_x & -\partial_{x^2}^2 & 0
    \end{bmatrix}}_{=: \mathcal{J}_r^E} \begin{bmatrix}
         \sigma_w \\ \sigma_\phi \\ v 
    \end{bmatrix}.
\end{equation}
    Note that by hypothesis, $\partial_t \varepsilon_{w,\phi} = \partial_t(\partial_x \phi - w) = 0$, hence $\partial_x\phi - w$ is constant over time, however, this constraint is not written explicitly in the previous set of equations.
\end{remark}

\subsection{Implicit representations}
\label{ss-ImpRep}
\subsubsection{Timoshenko}  Let us now write the Timoshenko beam equation in implicit form, i.e. by putting the differential operators in the Hamiltonian.
\begin{equation} \label{eqn:dynamic-timoshenko-implicit}
	\begin{bmatrix}
	\partial_t w \\
    \partial_t p_w \\
    \partial_t \phi \\
	\partial_t p_\phi 
	\end{bmatrix} = \begin{bmatrix}
	0 & \Id & 0 & 0 \\
	-\Id & 0 & 0 & 0 \\
	0 & 0 & 0 & \Id \\
	0 & 0 & -\Id & 0
\end{bmatrix}
\begin{bmatrix}
	\delta_w \mathcal{H}^I \\
	\delta_{p_w} \mathcal{H}^I \\
	\delta_\phi \mathcal{H}^I  \\
	\delta_{p_\phi} \mathcal{H}^I 
\end{bmatrix},
\end{equation}
with Hamiltonian 
\begin{equation}\label{eqn:ham-timoshenko-implicit}\begin{aligned} 
\mathcal{H}^I :=    &\frac{1}{2} \int_\Omega \Big(\frac{1}{\rho A}p_w^2  + \frac{1}{\rho I}p_{\phi}^2 + T_0 (\partial_x w)^2 \\&+ EI(\partial_x \phi)^2 + A \kappa G (\partial_x w - \phi)^2 \Big) \d x. 
\end{aligned}\end{equation}
Notice that boundary terms will then appear in the constitutive relations! Moreover, this formulation allows for direct access to the displacement $w$  instead of its gradient $\varepsilon_w$. Let us now write the constitutive relations by using Corollary~\ref{corollary-variational-derivative}
$$ 
\begin{cases}
    \delta_w \mathcal{H}^I = - \partial_x (T_0 \partial_x w) - \partial_x(A \kappa G (\partial_x w- \phi)), \\
	\delta_{p_w} \mathcal{H}^I = \frac{p_w}{\rho A}, \\
	 \delta_\phi \mathcal{H}^I = -  \partial_x (EI \partial_x \phi ) - A \kappa G (\partial_x w - \phi), \\
	 \delta_{p_\phi} \mathcal{H}^I = \frac{p_\phi}{\rho I}.
\end{cases} 
$$
We can now define the structure matrix $\mathcal{J}$ and the Lagrange subspace operators $\mathcal{S}$ and $\mathcal{P}$ as
$$ 
\mathcal{J}^I =  \begin{bmatrix}
	0 & \Id & 0 &0 \\
	-\Id & 0 & 0 & 0\\
	0 & 0 & 0 & \Id \\
	0 & 0 & -\Id & 0
\end{bmatrix}, \quad  \mathcal{P}^I =  \begin{bmatrix}
\Id & 0 & 0 & 0\\
0 & \Id & 0  & 0\\
0 & 0 & \Id & 0 \\
0& 0 & 0 & \Id
\end{bmatrix},
$$
$$
\mathcal{S}^I =  \begin{bmatrix}
\mathcal{S}_{1,1}  & 0 & \partial_x (A \kappa G \cdot ) &  0\\
0 & \frac{1}{\rho A}  & 0  & 0\\
 - A \kappa G \partial_x  \cdot  & 0 &  \mathcal{S}_{3,3} & 0 \\
0 & 0 &0&   \frac{1}{\rho I}
\end{bmatrix}, 
$$
with $ \mathcal{S}^I_{1,1} := - \partial_x(T_0 \partial_x \cdot ) - \partial_x(A\kappa G(\partial_x \cdot ))  $ and $\mathcal{S}^I_{3,3} :=- \partial_x(EI\partial_x \cdot ) + A \kappa G\,\Id  $.

Then, writing $z^I = ( w,p_w,\phi, p_\phi)^\top$, and $e^I = (\delta_w \mathcal{H},\delta_{p_w} \mathcal{H},  \delta_\phi \mathcal{H}, \delta_{p_\phi} \mathcal{H})^\top$, we get
$$ 
\begin{cases}
	\partial_t z^I = \mathcal{J}^I e^I, \\
	\mathcal{P}^{I\dag} e^I = \mathcal{S}^{I\dag} z^I.
\end{cases} 
$$
Let us finally write the system using the image representation
$$ 
\begin{cases}
	\mathcal{P}^I \partial_t \xi^I = \mathcal{J}^I e^I, \\
	 e^I = \mathcal{S}^I  \xi^I.
\end{cases} 
$$
Notice that, $\mathcal{S}^\dag \mathcal{P} = \mathcal{P}^\dag \mathcal{S} \geq 0$ which comes from the non-negativity of  the Hamiltonian \eqref{eqn:ham-timoshenko-implicit}.  Moreover, defining $\mathcal{R} = \begin{bmatrix}
    \mathcal{P} & \mathcal{S}
\end{bmatrix}^\top$, since $\mathcal{S}$ is invertible, the associated polynomial matrix of the operator $R(s)$ verifies $\text{rank}(R(s)) = 4 \, \, \forall s \in \mathbb{C}$ which gives the maximal reciprocity of the operator $\mathcal{R}$ (see \cite{maschke2023linear} for details). 

\begin{theorem}(Timoshenko power balance - implicit case) The power balance of \eqref{eqn:ham-timoshenko-implicit} reads
\begin{equation}
    \begin{aligned}
    \frac{\d }{\d t} \mathcal{H}^I = 
        [\,& \partial_t w\,(T_0\,\partial_x w + A\kappa G\,(\partial_xw - \phi)) \\& + \partial_t \phi\,EI \partial_x \phi \,]_a^b.
    \end{aligned}
\end{equation}
\end{theorem}
One can identify the \emph{\textbf{energy} boundary port} variable 
$$
\begin{aligned}\chi_\partial &= \text{tr}(\begin{bmatrix}
    w, && w, && \phi 
\end{bmatrix}^\top), \\ \varepsilon_\partial &= \text{tr}(\begin{bmatrix}
    T_0\,\partial_x w, && A\kappa G\,(\partial_x w - \phi), && EI\,\partial_x \phi
\end{bmatrix}^\top),\end{aligned}
$$ 
this yields $\frac{\d}{\d t} \mathcal{H}^I = [\,\frac{\d}{\d t} (\chi_\partial) \cdot \varepsilon_\partial\,]_a^b$.

\subsubsection{Euler--Bernoulli}



Let us reduce the previously defined system by putting $\frac{1}{\rho A}p_\phi^2 =0$ and $A \kappa G (\partial_x w - \phi)^2 = 0$ in the Hamiltonian \eqref{eqn:ham-timoshenko-implicit}. One gets: $ \partial_t p_\phi = 0, \quad \phi = \partial_x w$. This allows us to define the constrained Hamiltonian as
$$
\mathcal{H}^I_c = \frac{1}{2} \int_\Omega \left(  \frac{1}{\rho A}p_w^2  + T_0 (\partial_x w)^2 + EI(\partial_{x^2}^2 w)^2 \right) \d x.
$$

The constrained dynamic reads
\begin{equation} \label{eqn:dynamic-eulerbernoulli-implicit-DAE}
\begin{bmatrix}
    \Id &0 & 0  & 0\\
    0& \Id & 0 & 0\\
     \partial_x &0 & 0 & 0 \\
    0 & 0 & 0 &0\\
\end{bmatrix} \!\!
	\begin{bmatrix}
	\partial_t w \\
    \partial_t p_w \\
    \partial_t \phi \\
	\partial_t p_\phi 
	\end{bmatrix} \!=\! \begin{bmatrix}
	0 & \Id & 0 & 0 \\
	-\Id & 0 & 0 & 0 \\
	0 & 0 & 0 & \Id \\
	0 & 0 & -\Id & 0
\end{bmatrix}
\begin{bmatrix}
	\delta_w \mathcal{H}^I \\
	\delta_{p_w} \mathcal{H}^I \\
	\delta_\phi \mathcal{H}^I \\
	\delta_{p_\phi} \mathcal{H}^I 
\end{bmatrix}\!.
\end{equation}

\begin{theorem}(Euler--Bernoulli power balance - Implicit case)
\begin{equation}
\begin{aligned}
    \frac{\d }{\d t} \mathcal{H}^I_c = [\,& \partial_t w\,(T_0 \partial_x w - \partial_x(EI \partial_{x^2}^2 w)) \\ & + \partial_t(\partial_x w)\,EI\,\partial_{x^2}^2w   \,]_a^b.
\end{aligned}    
\end{equation}
\end{theorem}
One can identify the \emph{\textbf{energy} boundary port} variable 
$$
\begin{aligned}\chi_\partial^c &= \text{tr}(\begin{bmatrix}
    w, && w, && \partial_x w 
\end{bmatrix}^\top), \\ \varepsilon_\partial^c &= \text{tr}(\begin{bmatrix}
    T_0\,\partial_x w, && - \partial_x (EI\,\partial_{x^2}^2 w ), && EI\,\partial_{x^2}^2 w
\end{bmatrix}^\top),\end{aligned}
$$ 
this yields $\frac{\d}{\d t} \mathcal{H}^I = [\,\frac{\d}{\d t}(\chi_\partial^c) \cdot \varepsilon_\partial^c\,]_a^b$.

\begin{remark}
Solving the two constraints $\partial_t p_\phi$ and $\partial_x w = \phi$ yields the reduced unconstrained system
\begin{equation} \label{eqn:dynamic-eulerbernoulli-implicit}
	\begin{bmatrix}
		\partial_t w \\
		\partial_t p_w \\
	\end{bmatrix} = \underbrace{\begin{bmatrix}
	0 & \Id   \\
	-\Id & 0  \\
\end{bmatrix}}_{=:\mathcal{J}_r^I}
\begin{bmatrix}
	\delta_w \mathcal{H}^I_c \\
	\delta_{p_w} \mathcal{H}^I_c \\
\end{bmatrix},
\end{equation}
together with the constitutive relations
\begin{equation*}
    \begin{cases}
        \delta_w \mathcal{H}^I_c = - \partial_x (T_0 \partial_x w) + \partial_{x^2}^2(EI \partial_{x^2}^2 w). \\
        \delta_{p_w} \mathcal{H}^I_c = \frac{p_w}{\rho A}, \\
    \end{cases}
\end{equation*}
\end{remark}

\subsection{From implicit to explicit formulations}
\label{ss-EquivRep}
\subsubsection{Timoshenko}
In this subsection, we will denote by a $z^E,e^E$ the state and effort variables in the explicit formulation and by a $z^I,e^I$ the state and effort in the implicit formulation. 
${\cal J}^E,{\cal S}^E,{\cal P}^E$ denote the structure and constitutive matrices in the explicit case and ${\cal J}^I,{\cal S}^I,{\cal P}^I$ the structure and constitutive matrices in the implicit case. 
Let us now describe a procedure that allows to pass from  the implicit to the explicit representation, this procedure is similar to the one described in \cite[Section 6.]{mehrmann2023differential} but uses differential operators instead of matrices. 
First let us define a transformation operator $\mathcal{G}$ from $(H^1 \cap L^2_0) \times L^2 \times  H^1 \times  L^2$ to $(L^2)^5$, where $L^2_0$ is the space of square-integrable functions with zero mean,  $\int_\Omega w = 0$, as
\begin{equation*}
\mathcal{G} := \begin{bmatrix}
\partial_x  & 0 & 0 & 0 \\
0 & \Id & 0 & 0 \\
0 & 0 & \partial_x  & 0 \\
0 & 0 & 0 & \Id \\
\partial_x & 0 & -\Id & 0
\end{bmatrix}.
\end{equation*}
Its inverse is actually well defined, thanks to  Poincaré-Wirtinger's inequality, since $\int_\Omega w = 0$. Physically speaking, it means that one neglects rigid body motions of the beam. The inverse $\mathcal{F}$ of $\mathcal{G}$ is given by
\begin{equation}
    \mathcal{F}(z) := z
    \mapsto \!\!\!
    \begin{bmatrix}
        x \mapsto \left( \int_a^x z_2 \d y - \frac{1}{b-a} \int_a^b z_2 \d y \right) \\
        z_1 \\
        z_2 - z_5 \\
        z_3 \\
    \end{bmatrix}.
\end{equation}
A direct computation then yields
$$ 
\begin{cases}
	\mathcal{G} z^I = z^E, & 
	e^I =  \mathcal{G}^\dag  e^E, \\
    z^I  = \mathcal{F} z^E, &
    \mathcal{F}^\dag e^I = e^E.
\end{cases} 
$$
Note that a \emph{geometric} interpretation might be useful to understand these transformations since states are usually described as vectors and co-states  as covectors, hence they are contravariant and covariant, respectively.
Then, one can compute the time derivative in the implicit case and compare it to the explicit case
$$
\partial_t \mathcal{G}(z^I) = \mathcal{G}\mathcal{J}^Ie^I = \mathcal{G}\mathcal{J}^I\mathcal{G}^\dag e^E.
$$
A direct computation shows that $\mathcal{G}\mathcal{J}^I\mathcal{G}^\dag = \mathcal{J}^E$.
Additionally, regarding the Lagrange subspace
$$ 
\mathcal{S}^{E\top} z^E = \mathcal{P}^{E\top} e^E.
$$
Let us premultiply this equation by $\mathcal{G}^\dag$ and replace the state and effort variables by their explicit counterparts. We get to
$$
\mathcal{G}^\dag \mathcal{S}^{E\top} \mathcal{G} z^I = \mathcal{G}^\dag \mathcal{P}^{E\top} \mathcal{F}^\dag e^I,
$$
and we have $\mathcal{G}^\dag \mathcal{S}^{E\top} \mathcal{G} = \mathcal{S}^{I\dag}$ and $\mathcal{G}^\dag \mathcal{P}^{E\top} \mathcal{F}^\dag = \mathcal{P}^{I\dag}$.

\subsubsection{From Timoshenko to Euler--Bernoulli}

\begin{figure}[ht]
    \centering
    \begin{tikzpicture}
        \node (A) at (-2,1) {\eqref{eqn:dynamic-timoshenko-explicit}, $\mathcal{H}^E$};
        \node (B) at (2,1) {\eqref{eqn:dynamic-timoshenko-implicit}, $\mathcal{H}^I$};
        \node (C) at (-2,-0.5) {\eqref{eqn:dynamic-eulerbernoulli-explicit-DAE}, $\mathcal{H}_c^E$};
        \node (D) at (2,-0.5) {\eqref{eqn:dynamic-eulerbernoulli-implicit-DAE}, $\mathcal{H}_c^I$};

        \draw[<->] (A) -- (B)  node[midway,above]{$\mathcal{G},\mathcal{F}$};
        
        \draw[->] (A) -- (C)  node[midway,left]{$\Pi^E$};
        
        \draw[->] (B) -- (D)  node[midway,right]{$\Pi^I$};
        
        \draw[<->] (C) -- (D)  node[midway,below]{$\mathcal{G},\mathcal{F}$};
    \end{tikzpicture}
    \caption{Diagram of the transformations between beam models and representations}
    \label{fig:diagram-transformations}
\end{figure}
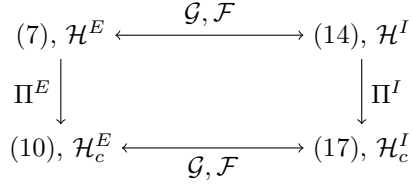

Let us denote
\begin{equation}
    \Pi^E := \begin{bmatrix}
        \Id & 0 & 0 & 0 & 0 \\
        0 & \Id & 0 & 0 & 0 \\
        0 & 0 & \Id & 0 & 0 \\
        0 & 0 & 0 & 0 & 0 \\
        0 & 0 & 0 & 0 & 0
    \end{bmatrix},
\end{equation}
the (singular) operator that allows to constrain the system \eqref{eqn:dynamic-timoshenko-explicit}; when premultiplying the left side of equation \eqref{eqn:dynamic-timoshenko-explicit} with $\Pi^E$ one gets to the constrained system \eqref{eqn:dynamic-eulerbernoulli-explicit-DAE}. In order to get the implicit counterpart $\Pi^I$, let us apply $\mathcal{G}$ to pass from \eqref{eqn:dynamic-timoshenko-implicit} to \eqref{eqn:dynamic-timoshenko-explicit}, then $\Pi^E$ to pass from \eqref{eqn:dynamic-timoshenko-explicit} to \eqref{eqn:dynamic-eulerbernoulli-explicit-DAE}, and finally $\mathcal{F}$ to pass from \eqref{eqn:dynamic-eulerbernoulli-explicit-DAE} to \eqref{eqn:dynamic-eulerbernoulli-implicit-DAE}. Therefore, let us compute $\mathcal{F} \, \Pi^E \, \mathcal{G}$
$$
\Pi^I := \mathcal{F} \, \Pi^E \, \mathcal{G} = \begin{bmatrix}
    \Id & 0 & 0 & 0\\
    0 & \Id & 0 & 0\\
    \partial_x & 0 & 0 & 0 \\
    0 & 0 & 0 & 0 
\end{bmatrix}.
$$
Finally, one can represent these transformations as a commutative diagram, shown in Fig. \ref{fig:diagram-transformations}.

\subsubsection{Euler--Bernoulli}

Let us now apply  the same procedure to the Euler--Bernoulli \emph{reduced} case, i.e., let us pass from \eqref{eqn:dynamic-eulerbernoulli-explicit} to \eqref{eqn:dynamic-eulerbernoulli-implicit}. Let us denote by 
$$
\mathcal{G}_r := \begin{bmatrix}
    \Id &0 &0  \\
     0& \partial_x  & \partial_{x^2}^2 \\
\end{bmatrix}^\top,
$$ 
the transformation operator between the two reduced systems.

Then, assuming $\varepsilon_\phi = \partial_x w$, i.e., on the subspace where this constraint is satisfied, we have
$$
\mathcal{G}_r \begin{bmatrix}
    w \\ p_w
\end{bmatrix} = \begin{bmatrix}\varepsilon_w \\ \varepsilon_\phi \\
    p_w 
\end{bmatrix},
\quad
\text{and}
\quad
\mathcal{G}_r^\dag \begin{bmatrix}
    \sigma_w \\
    \sigma_\phi \\
    v 
\end{bmatrix} = \begin{bmatrix}
    \delta_w \mathcal{H}^I \\
    \delta_{p_w} \mathcal{H}^I
\end{bmatrix}.
$$
Additionally, one has: $
\mathcal{G}_r\mathcal{J}_r^I \mathcal{G}_r^\dag = \mathcal{J}_r^E.
$



\section{Conclusion \& Outlook}

Implicit Stokes-Lagrange representations have been proposed for two examples of pHs, respectively with local and non local dissipation. The corresponding implicit Hamiltonian and power boundary and energy boundary port-variables have been derived. \\
Stokes-Lagrange representations have also been derived for the Timoshenko and the Euler-Bernoulli models. Two bijective transformations between the corresponding explicit and implicit models have been exhibited. Since the Euler-Bernoulli beam is a flow-constrained version of the Timoshenko beam, it has been checked that the corresponding projection operators commute with the transformations between explicit and implicit representations of both models. 

We will now consider extending these results to 2D (Dzekster seapage model, Reissner-Mindlin and Kirchoff-Love plate models) and 3D (Maxwell equations) examples. Another avenue we want to explore concerns the comparative analysis of explicit and implicit formulations, in terms of their numerical properties.

\bibliographystyle{plain}
\bibliography{mtns2024}

\begin{thebibliography}{10}

\bibitem{beattie2018linear}
Christopher Beattie, Volker Mehrmann, Hongguo Xu, and Hans Zwart.
\newblock Linear port-{H}amiltonian descriptor systems.
\newblock {\em Mathematics of Control, Signals, and Systems}, 30:1--27, 2018.

\bibitem{bendimerad2023implicit}
Antoine Bendimerad-Hohl, Ghislain Haine, Laurent Lef{\`e}vre, and Denis
  Matignon.
\newblock Implicit port-{H}amiltonian systems: structure-preserving
  discretization for the nonlocal vibrations in a viscoelastic nanorod, and for
  a seepage model.
\newblock {\em IFAC-PapersOnLine}, 56(2):6789--6795, 2023.

\bibitem{BHMCAM2023}
Andrea Brugnoli, Ghislain Haine, and Denis Matignon.
\newblock {Stokes-Dirac structures for distributed parameter port-Hamiltonian
  systems: An analytical viewpoint}.
\newblock {\em Communications in Analysis and Mechanics}, 15(3):362--387, 2023.

\bibitem{ducceschi2019conservative}
Michele Ducceschi and Stefan Bilbao.
\newblock Conservative finite difference time domain schemes for the
  prestressed {T}imoshenko, shear and {Euler--B}ernoulli beam equations.
\newblock {\em Wave Motion}, 89:142--165, 2019.

\bibitem{duindam2009modeling}
Vincent Duindam, Alessandro Macchelli, Stefano Stramigioli, and Herman
  Bruyninckx.
\newblock {\em Modeling and control of complex physical systems: the
  port-{H}amiltonian approach}.
\newblock Springer Science \& Business Media, 2009.

\bibitem{Dzektser72}
E.~S. Dzektser.
\newblock Generalization of the equation of motion of ground waters with free
  surface.
\newblock {\em Dokl. Akad. Nauk SSSR}, 202(5):1031--1033, 1972.

\bibitem{heidari2019port}
Hanif Heidari and H~Zwart.
\newblock Port-{H}amiltonian modelling of nonlocal longitudinal vibrations in a
  viscoelastic nanorod.
\newblock {\em Mathematical and computer modelling of dynamical systems},
  25(5):447--462, 2019.

\bibitem{heidari2022nonlocal}
Hanif Heidari and Hans Zwart.
\newblock Nonlocal longitudinal vibration in a nanorod, a system theoretic
  analysis.
\newblock {\em Mathematical Modelling of Natural Phenomena}, 17:24, 2022.

\bibitem{jacob2022solvability}
Birgit Jacob and Kirsten Morris.
\newblock On solvability of dissipative partial differential-algebraic
  equations.
\newblock {\em IEEE Control Systems Lett.}, 6:3188--3193, 2022.

\bibitem{jacob2012linear}
Birgit Jacob and Hans~J Zwart.
\newblock {\em Linear port-{H}amiltonian systems on infinite-dimensional
  spaces}, volume 223.
\newblock Springer Science, 2012.

\bibitem{maschke2023linear}
Bernhard Maschke and Arjan van~der Schaft.
\newblock Linear boundary port-{H}amiltonian systems with implicitly defined
  energy.
\newblock {\em arXiv preprint arXiv:2305.13772}, 2023.

\bibitem{mehrmann2023control}
Volker Mehrmann and Benjamin Unger.
\newblock Control of port-{H}amiltonian differential-algebraic systems and
  applications.
\newblock {\em Acta Numerica}, 32:395--515, 2023.

\bibitem{mehrmann2023differential}
Volker Mehrmann and Arjan van~der Schaft.
\newblock Differential--algebraic systems with dissipative {H}amiltonian
  structure.
\newblock {\em Mathematics of Control, Signals, and Systems}, pages 1--44,
  2023.

\bibitem{preuster2024jet}
Till Preuster, Manuel Schaller, and Bernhard Maschke.
\newblock Jet space extensions of infinite-dimensional {H}amiltonian systems.
\newblock {\em arXiv preprint arXiv:2401.15096}, 2024.

\bibitem{rashad2020twenty}
Ramy Rashad, Federico Califano, Arjan~J van~der Schaft, and Stefano
  Stramigioli.
\newblock Twenty years of distributed port-{H}amiltonian systems: a literature
  review.
\newblock {\em IMA Journal of Mathematical Control and Information},
  37(4):1400--1422, 2020.

\bibitem{SchoberlSchlacherMathMod2015}
M.~Schöberl and K.~Schlacher.
\newblock {Lagrangian and port-Hamiltonian formulation for
  distributed-parameter systems}.
\newblock {\em IFAC-PapersOnLine}, 48(1):610--615, 2015.

\bibitem{scholberl2014auto}
Markus Schöberl and Andreas Siuka.
\newblock Jet bundle formulation of infinite-dimensional port-{H}amiltonian
  systems using differential operators.
\newblock {\em Automatica}, 50:607--613, 2 2014.

\bibitem{van2014port}
Arjan van~der Schaft, Dimitri Jeltsema, et~al.
\newblock Port-{H}amiltonian systems theory: An introductory overview.
\newblock {\em Foundations and Trends{\textregistered} in Systems and Control},
  1(2-3):173--378, 2014.

\bibitem{van2018generalized}
Arjan van~der Schaft and Bernhard Maschke.
\newblock Generalized port-{H}amiltonian {DAE} systems.
\newblock {\em Systems \& Control Letters}, 121:31--37, 2018.

\bibitem{van2020dirac}
Arjan van~der Schaft and Bernhard Maschke.
\newblock Dirac and {L}agrange algebraic constraints in nonlinear
  port-{H}amiltonian systems.
\newblock {\em Vietnam Journal of Mathematics}, 48:929--939, 2020.

\bibitem{van2002hamiltonian}
Arjan~J van~der Schaft and Bernhard~M Maschke.
\newblock {H}amiltonian formulation of distributed-parameter systems with
  boundary energy flow.
\newblock {\em Journal of Geometry and Physics}, 42(1-2):166--194, 2002.

\bibitem{yaghi2022port}
Mohammed Yaghi, Fran{\c{c}}oise Couenne, Aur{\'e}lie Galfr{\'e}, Laurent
  Lef{\`e}vre, and Bernhard Maschke.
\newblock Port-{H}amiltonian formulation of the solidification process for a
  pure substance: A phase field approach.
\newblock {\em IFAC-PapersOnLine}, 55(18):93--98, 2022.

\end{thebibliography}

\end{document}